\documentclass[12pt, 14paper,reqno]{amsart}
\usepackage{amsmath,amsfonts,amssymb}
\usepackage{longtable}

\theoremstyle{plain}
\newtheorem{theorem}{Theorem}
\newtheorem{lemma}{Lemma}
\newtheorem{corollary}{Corollary}[section]
\newtheorem{proposition}{Proposition}[section]
\theoremstyle{proof}
\theoremstyle{definition}

\theoremstyle{remark}
\newtheorem{remark}{Remark}
\theoremstyle{lamma}

\numberwithin{equation}{section}
\numberwithin{lemma}{section}
\numberwithin{theorem}{section}

\usepackage{amssymb, amsmath, amsthm}
\usepackage[breaklinks]{hyperref}
\theoremstyle{thmrm}

\usepackage{graphicx}

\begin{document} 
\title[Generalized Lebesgue-Ramanujan-Nagell equations]{On the solutions of certain  Lebesgue-Ramanujan-Nagell  equations}
\author{Kalyan Chakraborty, Azizul Hoque and Richa Sharma}
\address{Kalyan Chakraborty @Kerala School of Mathematics, Kozhikode 673571, Kerala, India.}
\email{kalychak@ksom.res.in}
\address{Azizul Hoque @Department of Mathematics, Rangapara College, Rangapara, Sonitpur 784505, Assam, India.}
\email{ ahoque.ms@gmail.com}
\address{Richa Sharma @Malaviya National Institute of Technology Jaipur,  Jaipur 302017, India.}
\email{richasharma582@gmail.com}
\keywords{Diophantine equation, Integer solution, S-integers, Lucas sequences, Primitive divisors.}
\subjclass[2010] {Primary: 11D41, 11D61, Secondary: 11R29.}
\maketitle
     \begin{abstract}
     We completely solve the Diophantine equation $x^2+2^k11^\ell19^m=y^n$ in integers $x,y\geq 1;~ k,\ell, m\geq 0~$ and $n\geq 3$ with $\gcd(x,y)=1$, except the case $2\mid k, 2\nmid \ell m$ and $5\mid n$. We use this result to  recover some earlier results in the same direction. 
     \end{abstract}
     \maketitle

\section{Introduction}\noindent 
Generalized Lebesgue-Ramanujan-Nagell equation
\begin{equation}\label{eqi}
x^2+C=\lambda y^n, ~x,y, \lambda\geq 1, n\geq 3, 
\end{equation}
with $C$ a fixed positive integer is being studied by many authors.
The earliest result is due to Fermat who  proved that $(x,y)=(5,3)$ is the only solution of \eqref{eqi} when $(C, \lambda,n)=(2,1,3)$.  Lebesgue \cite{LE1850} in 1850 proved that \eqref{eqi} has no solution when $(C, \lambda)=(1,1)$. In 1993, Cohn \cite{CO93} investigated \eqref{eqi} for $\lambda=1$, and he solved the equation for $77$ values of $C$ in the range $1\leq C\leq 100$. The remaining values of $C$ in this range were covered in \cite{BM20, MW96} (when $\lambda=1$). There are numerous results for $\lambda=2, 4$, and interested readers can look into \cite{BHS18, CHS-21, LTT, MFST} and references therein. In the recent years, several authors became interested in the case when $C = p_1^{a_1}p_2^{a_2}\cdots p_k^{a_k}$, where $p_i$'s are distinct primes and $k\geq 1, ~ a_i\geq 0$ are integers. While  $(\lambda,C)=(1,p_1^{a_1})$ the equation \eqref{eqi} is completely settled by the combined work of several authors (cf. \cite{AM02, BP08, BP, CO92, HS16, LE02}).  For $(\lambda, C)=(2,p_1^{a_1})$ with even $a_1$, \eqref{eqi} is studied in \cite{HT20, TE}. Recently, Hoque studied a more general Lebesgue-Ramanujan-Nagell equation in \cite{HO20}. On the other hand, Chakraborty et al. \cite{CHS-19} studied \eqref{eqi} when $(\lambda, C)=(4, p_1^{a_1})$, and completely  solved it for some fixed values of $p_1$. There are some results concerning the solution of \eqref{eqi} when $\lambda=1$ and $k=2$, i.e. when  $C=p_1^{a_1}p_2^{a_2}$ for some primes $p_1, p_2\leq 19$. For instance, all the solutions of \eqref{eqi} when $\gcd(x, y)=1$ and $C=2^a5^b, 2^a13^b, 5^a13^b, 2^a19^b$ were obtained in \cite{ MLT08, LT08, LT09, SUZ12}.  A natural next step in this series of explorations is to consider \eqref{eqi} when $C$ has three distinct prime factors. Recent progress in this case were made in \cite{CDILS13, GH19, GDT16, GLT08} for $C=2^a3^b11^c, 2^a3^b13^c, 2^a3^b17^c, 2^a13^b17^c, 2^a5^b13^c$.  

Here we completely solve,
\begin{equation}\label{me}
x^2+2^k11^\ell 19^m=y^n
\end{equation}
in positive integers $x, y$  with $\gcd(x,y)=1$ and in non-negative integers $k,\ell, m, n$, except the case $2\mid k, 2\nmid \ell m$ and $5\nmid n$.
We use this result to deduce the results of Cang\"{u}l et. al. \cite{CDLPS10}, Soydan et. al. \cite{SUZ12}, and some new results. 

\section{Statement of results}
We begin with the case $k=\ell=m=0$, that is,
$$
x^2+1=y^n.
$$
It was proved in \cite{LE1850} that this equation has no solution in positive integers $x, y, n$, when $n\geq 3$.  Clearly, it has no positive integer solution for $n=2$.  When $n=1$, it has  a family of solutions $(x,y)=(t,t^2+1)$.
It is easy to see that \eqref{me} has no solution in positive integers $x, y, k, \ell, m $ when $n=0$.

The first non-trivial case is when $n=1$. In this case, \eqref{me}  is parametrized simply by 
\begin{eqnarray*}
(x,y)= \left(t, t^2+2^k 11^\ell 19^m\right).
\end{eqnarray*}

Moving to $n=2$ and  \eqref{me} becomes
\begin{equation}\label{me2} x^2+2^k11^\ell19^m=y^2.\end{equation}
Let $x$ be even. Then $y$ is odd (since $\gcd(x,y)=1$) and that implies $k=0$. Thus   
\eqref{me2} reduces to 
$$
(y+x)(y-x)=11^\ell 19^m.
$$
Now  $\gcd(x,y)=1$ would force $\gcd(y-x,y+x)=1$. Using this fact, we obtain the following solutions of this equation:
$$(x,y)=\left(\frac{11^{\ell-\ell_1}19^{m-m_1}-11^{\ell_1}19^{m_1}}{2}, \frac{11^{\ell-\ell_1}19^{m-m_1}+11^{\ell_1}19^{m_1}}{2}\right),$$
where $(\ell_1,m_1)\in\{0,\ell\}\times \{0, m\}\setminus\{(\ell, m)\}$ with $\ell\equiv m\pmod 2$.

If $x$ is considered odd then $y$ is even and as before $k=0$. Therefore the solutions are given as in the last case with the condition $\ell\not\equiv m\pmod 2$. 

Finally if both $x$ and $y$ are odd, then $\gcd(y-x,y+x)=2$ and that gives $k\geq 2$. For $k=2$,  \eqref{me2} can be written as 
$$
X^2+X+11^\ell19^m=Y^2+Y,
$$
where $x=2X+1$ and $y=2Y+1$. The impossibility of this can be easily obtained by reading it modulo $2$, and hence $k\geq 3$. Thus \eqref{me2} can be written as 
$$\left(\frac{y-x}{2}\right)\left(\frac{y+x}{2}\right)=2^{k-2}11^\ell19^m.$$
Since $\gcd\left(\frac{y-x}{2},\frac{y+x}{2}\right)=1$, so that we get
$$(x,y)=\left( 2^{k-k_1-2}11^{\ell-\ell_1}19^{m-m_1}-2^{k_1}11^{\ell_1}19^{m_1}, 2^{k-k_1-2}11^{\ell-\ell_1}19^{m-m_1}+2^{k_1}11^{\ell_1}19^{m_1}   \right),$$
where $(k_1,\ell_1,m_1)\in\{0,k-2\}\times\{0, \ell\}\times\{0, m\}\setminus\{(k-2, \ell, m)\}$.
 
In order to complete the study of \eqref{me}, it remains to consider the following equation:
\begin{equation}\label{men}
x^2+2^k11^\ell 19^m=y^n, ~x, y\geq 1, ~ \gcd(x,y)=1,~ k,\ell,m\geq 0 ~ n\geq 3.
\end{equation}
The main result of this paper is:
\begin{theorem}\label{thm} 
 \eqref{men} has no solutions $(x, y, k, \ell, m, n)$ with $n \not\in \{3, 4, 5, 6, 10, 12\}$ except when $2\mid k, 2\nmid \ell m, 5\mid n$. 
For $n \in \{3, 4, 5, 6, 10, 12\}$, all the solutions are given by Tables \ref{T3}, \ref{T612}, \ref{T4},  \ref{Tp} and $(x,y,k,\ell,m,n)=(241,3,3,2,0,10)$.
\end{theorem}
As it is with other such equations, one needs to consider case by case to prove this theorem. While passing through various  cases many a times we use MAGMA and in some cases PARI (v.2.9.1) to get to the desired results.

\section{Proof of Theorem \ref{thm}}
Assume that $y$ is even then $x$ is odd 
and thus \eqref{men} gives $k=0$. Therefore it follows that $x^2+11^\ell 19^m\equiv 2, 4\pmod 8$. However $y^n\equiv 0 \pmod 8$ since $n\geq 3$, and hence $y$ is odd. We first consider $n=3$, that is, 
\begin{equation}\label{me3}
x^2+2^k 11^\ell 19^m=y^3,~ x, y\geq 1, \gcd(x,y)=1, k, \ell, m\geq 0.
\end{equation}      

\begin{proposition}\label{3p}
All the solutions of \eqref{me3} are given in Table \ref{T3}.
\end{proposition}

\begin{center}
{\small
\begin{longtable}{ccccc | ccccc}
\caption{Solutions of \eqref{men} for $n=3$} \label{T3} \\
\hline \multicolumn{1}{c}{$x$} & \multicolumn{1}{c}{$y$} & \multicolumn{1}{c}{$k$}& \multicolumn{1}{c}{$\ell$} & \multicolumn{1}{c}{$m$}& \multicolumn{1}{|c}{$x$} & \multicolumn{1}{c}{$y$} & \multicolumn{1}{c}{$k$}& \multicolumn{1}{c}{$\ell$} & \multicolumn{1}{c}{$m$}\\ \hline 
\endfirsthead

\multicolumn{10}{c}%
{{\bfseries \tablename\ \thetable{} -- continued from previous page}} \\
\hline \multicolumn{1}{c}{$x$} & \multicolumn{1}{c}{$y$} & \multicolumn{1}{c}{$k$}& \multicolumn{1}{c}{$\ell$} & \multicolumn{1}{c}{$m$}& \multicolumn{1}{|c}{$x$} & \multicolumn{1}{c}{$y$} & \multicolumn{1}{c}{$k$}& \multicolumn{1}{c}{$\ell$} & \multicolumn{1}{c}{$m$}\\ \hline 
\endhead
\hline \multicolumn{10}{c}{{Continued on next page}} \\ \hline
\endfoot
\hline 
\endlastfoot
18 & 7  & 0 & 0 & 1 &   5 & 9  & 6 & 1 & 0\\ 
4 & 3  & 0 & 1 & 0 &   6179 & 345  & 18 & 1 & 0\\ 
58&15&0&1&0&2&5&0&2&0\\
835&89&6&2&0&404003&5465&12&2&0\\
9324&443&0&3&0&5&3&1&0&0\\
2345&177&7&0&2&1265&123&1&0&4\\
65&23&1&1&2&11290511&50337&11&6&1\\
487&81&7&2&1&87&23&1&2&1\\
70953&1841&19&2&1&1509&137&7&2&1\\
315&47&1&2&1&195499 & 3369&13&2&1\\
12653&543&1&2&1&56879&1479&1&2&1\\
202323621&344639&1&8&1&9158107626651&437727497&19&2&7\\
159&119&1&2&3&250173&3971&1&2&4\\
3922311&26129&13&2&5&1861&159&1&4&1\\
79153&1863&1&4&3&5497&785&8&6&0\\
11&5&2&0&0&7&17&8&0&1\\
7&5&2&0&1&9&5&2&1&0\\
81711&1885&2&4&2&147743&3393&9&6&1\\
69&17&3&0&1&12419&537&15&0&1\\
293493&56813&3&6&3&4103&273&9&0&3\\
1217&201&3&2&3&1349867&15273&11&2&3\\
1006153&10041&3&2&3&117&25&4&2&0\\
9959&465&10&3&0&11&9&5&0&1\\
80363& 3273&5&0&7&62129&1569&17&0&1\\
733326857&944241&23&4&5&45&17&3&0&2\\
1015&101&2&0&1&5805&323&1&2&0\\
\end{longtable}}
\end{center}
\vspace{-9mm}

\begin{proof}
Let us write $2^k11^\ell 19^m=Dz^6$. Then $D$ and $z$ are of the form $2^a11^b19^c$ with $0\leq a,b,c\leq 5$ and $2^\alpha11^\beta 19^\gamma$ with $0\leq \alpha\leq \frac{k}{6}, 0\leq \beta\leq \frac{\ell}{6}, 0\leq \gamma\leq \frac{m}{6}$ respectively.  Thus \eqref{me3} can be written as
\begin{equation}\label{3pe1}
X^2=Y^3-D
\end{equation}
where $X=x/z^3, ~ Y=y/z^2$ and $D$ is sixth-power free. Assume that $S =\{2,11,19\}$. 
In order to find the integer solutions of \eqref{me3}, one needs to determine all $S$-integral points on 216 elliptic curves defined by  \eqref{3pe1}. We use MAGMA to determine these $S$-integral points. We obtain $54$ such $S$-integral points which provide desired solutions of \eqref{me3} and all of them are listed in Table \ref{T3}. Note that, we also obtain $7$ (resp.  $54$) $S$-integral points with $x=0$ (resp.  $\gcd(x,y)>1)$, that are not mentioned here.     
\end{proof}
We now consider $n$'s which are multiple of $3$. In this case,  we derive the following result from Table \ref{T3}. 
\begin{corollary}\label{cor3.1}
If $n> 3$ is a multiple of $3$ 
and  \eqref{men} has an integer solution $(x,y,k,\ell, m,n)$, then $n = 6, 12$. For these values, all the solutions of \eqref{men} are given in Table \ref{T612}.
\end{corollary}
\begin{table}[ht]
 \centering
\begin{tabular}{ c c c c c c |  c c c c c c} 
 $x$ & $y$ & $k$ & $\ell$ & $m$ & $n$ & $x$ & $y$ & $k$ & $\ell$ & $m$ & $n$\\
\hline 
5&3&6&1&0&6&
487&9&7&2&1&6\\
117&5&4&2&0&6&
487&3&7&2&1&12\\
11&3&5&0&1&6\\
\hline
\end{tabular} \vspace{1mm}
\caption{\small Solutions of \eqref{men} when $n=6,12.$} \label{T612}
\end{table}

We now consider \eqref{men} when $n=4$, that is, 
\begin{equation}\label{me4}
x^2+2^k 11^\ell 19^m=y^4,~ x, y\geq 1, \gcd(x,y)=1, k, \ell, m\geq 0.
\end{equation}      
\begin{proposition}\label{4p}
All the solutions of \eqref{me4} are given in Table \ref{T4}.
\end{proposition}
\begin{table}[ht]
 \centering
\begin{tabular}{ c c c c c |  c c c c c} 
 $x$ & $y$ & $k$ & $\ell$ & $m$ & $x$ & $y$ & $k$ & $\ell$ & $m$ \\
\hline 
7&3&5&0&0&193&15&6&1&1\\
27&7&3&1&1&487&27&7&2&1\\
1293&37&3&3&1&&&&&\\
\hline
\end{tabular} \vspace{1mm}
\caption{\small Solutions of \eqref{me4}.} \label{T4}
\end{table}
\vspace{-9mm}
\begin{proof}
Assume that $2^k11^\ell 19^m=Dz^4$. Then $D$ and $z$ are of the form $2^a11^b19^c$ with $0\leq a,b,c\leq 3$ and $2^\alpha11^\beta 19^\gamma$ with $0\leq \alpha\leq \frac{k}{4}, 0\leq \beta\leq \frac{\ell}{4}, 0\leq \gamma\leq \frac{m}{4}$ respectively. Substituting $X=x/z^2$ and $ Y=y/z$ in \eqref{me4}, we get
\begin{equation}\label{4pe1}
X^2=Y^4-D,
\end{equation}
where  $D$ is fourth power-free. Since all the prime factors of $z$ are in $S=\{2,11,19\}$, so that
in order to find the solutions of \eqref{me4} one needs to determine all $S$-integral points on 64 elliptic curves defined by  \eqref{4pe1} and  MAGMA helps to determine these points. We obtain $4$ such
$S$-integral 
points which provide desired solutions of \eqref{me4} and they are listed in Table \ref{T4}. Note that we also obtain $1$ (resp. $4$) $S$-integral points which give $x=0$ (resp. $\gcd(x,y)>1)$.        
\end{proof}
\begin{remark}
Note that if $n \geq 4$ is a multiple of $4$ 
and \eqref{men} has an integer solution $(x,y,k,\ell, m,n)$, then $n=12$. This case is already considered in Corollary \ref{cor3.1}.
\end{remark}

\begin{proposition}\label{mp}
Let $p\geq 5$ be a prime. Then all the solutions of 
\begin{equation}\label{pe1}
x^2+2^k 11^\ell 19^m=y^p, ~x, y\geq 1,~\gcd(x,y)=1,  k, \ell, m\geq 0
\end{equation}
are given in Table \ref{Tp}, except for the case $2\mid k, 2\nmid \ell m$ with $p\ne 5$.
\end{proposition}

\begin{table}[ht]
 \centering
\begin{tabular}{ c c c c c c | c c c c c c} 
$x$ & $y$ & $k$ & $\ell$ & $m$ & $p$ & $x$ & $y$ & $k$ & $\ell$ & $m$ & $p$\\
\hline 
1&3&1&2&0&5& 22434&55&0&0&1&5\\
41&5&2&0&2&5& 241&9&3&2&0&5\\
\hline
\end{tabular} \vspace{1mm}
\caption{\small Solutions of \eqref{pe1}.} \label{Tp}
\end{table}

 We need few lemmas to prove Proposition \ref{mp}. It is easy to observe from the equations \eqref{pe5}, \eqref{pe7} and \eqref{pe7} that $1\leq x\leq 18, ~ 0\leq k\leq 7, ~ 0\leq \ell\leq 2, 0\leq m\leq 1$; $1\leq x\leq 22434, ~ 0\leq k\leq 28, ~ 0\leq \ell\leq 8, 0\leq m\leq 6$ and  $1\leq x\leq 279, ~ 0\leq k\leq 16, ~ 0\leq \ell\leq 4, 0\leq m\leq 3$, respectively. Using these bounds, one can get the following lemmas.

\begin{lemma}\label{pl1}
$(x,k,\ell, m)=(1,1,2,0)$ is the only solution of 
\begin{equation}\label{pe5}
x^2+2^k11^\ell19^m=3^5.
\end{equation}   
\end{lemma}
\begin{lemma}\label{pl2}
The equation
\begin{equation}\label{pe6}
x^2+2^k11^\ell19^m=55^5
\end{equation} 
has no solution, except $(x,k,\ell, m)=(22434, 0, 0, 1)$.  
\end{lemma}
\begin{lemma}\label{pl3}
The  equation
\begin{equation}\label{pe7}
x^2+2^k11^\ell19^m=5^7
\end{equation} 
has no solution.  
\end{lemma}
We are now in a position to give the proof of Proposition \ref{mp}.
\begin{proof}[\bf{Proof of  Proposition \ref{mp}}]
Let's recall that $y$ is odd.  The first step is to write 
 \eqref{pe1} as 
\begin{equation}\label{pe2}
x^2+dz^2=y^p,
\end{equation}
where $d \in \{ 1, 2,11,19,22,38,209,418\} $ depending on the exponents $k, \ell, m$, and $z=2^a11^b19^c$ for some non-negative integers $a, b$ and  $c$. 

Let $F=\mathbb{Q}(\sqrt{-d})$ and $\mathcal{O}_F$ its ring of integers. Next step is to factorize \eqref{pe2} in $\mathcal{O}_F$ as: 
$$
(x+z\sqrt{-d})(x-z\sqrt{-d})=y^p.
$$
Let's conclude that the ideals $(x+z\sqrt{-d})$ and $(x-z\sqrt{-d})$ are coprime in $\mathcal{O}_F$. If a prime ideal $\mathfrak{p}\subset \mathcal{O}_F$ is a common divisor of the ideals 
 then $x\pm z\sqrt{-d}\in \mathfrak{p}$ and thus $y\in \mathfrak{p}$.  Further $2x\in \mathfrak{p}$  implies either $2\in \mathfrak{p} $ or $x \in \mathfrak{p}$. Therefore $1$ lies in the ideals $<2, y>$ and $<x,y>$ since $\gcd(2, y)=\gcd(x,y)=1$. This shows that $1\in \mathfrak{p}$ which is not possible.  
Therefore by the unique factorization of ideals, $(x+z\sqrt{-d})$ can be expressed as a $p$-th power of some ideal $\mathfrak{I}\subset{O}_F$. On the other hand, $2, 3$ and $5$ are the only prime divisors of the class number $h(d)$ of $F$. For the exceptional case $2\mid k, 2\nmid \ell m$, we have $d=205$ and thus $h(-d)=20$. However in this case we assume that $p\ne 5$. Therefore $\gcd(p, h(d))=1$ and so  $\mathfrak{I}$ is a principal ideal. Therefore 
$$
x+z\sqrt{-d}=\varepsilon \alpha^p,
$$ 
for some $\alpha \in \mathcal{O}_F$ and $\varepsilon$ is a unit in $F$.
Further $\varepsilon$ can be absorbed into $\alpha^p$ as the order of the unit group is coprime to $p$. Thus we have 
\begin{equation} \label{pe3}
x+z\sqrt{-d}=\alpha^p.
\end{equation} 
It is well known that 
$$
\alpha=\begin{cases}
\frac{u+v\sqrt{-d}}{2}, ~~ u\equiv v\equiv 1\pmod 2 & \text{ if } d=11, 19,\\
u+v\sqrt{-d} & \text{ otherwise}, 
\end{cases}
$$ 
where both $u$ and $v$ are integers.
Thus one gets
\begin{equation}\label{pe4}
\frac{\alpha^p-\bar{\alpha}^p}{\alpha-\bar{\alpha}}=\begin{cases}
\frac{2z}{v} & \text{ if } d=11, 19,\\
\frac{z}{v} & \text{ otherwise}.
\end{cases}
\end{equation}

Let us define the Lucas sequence $\big(\mathcal{L}_n\big)_{n\geq 0}$ of integers by  
$$
\mathcal{L}_n=\frac{\alpha^n-\bar{\alpha^n}}{\alpha-\bar{\alpha}}, \text{ for } n\geq 0.
$$ 
 Let $P(t)$(for a non-zero integer $t$)
 be the largest prime factor of $t$ with the convention that $P(\pm1)=1$. Then by \eqref{pe4}
$P(\mathcal{L}_p) \leq 19$.
Let us recall that a prime divisor of $\mathcal{L}_n$ is called primitive if it does not divide $ (\alpha-\bar{\alpha})^2\mathcal{L}_1\mathcal{L}_2\cdots\mathcal{L}_{n-1}$. One of the most important property of a primitive divisor $q$ of $\mathcal{L}_n$ is that $q\equiv 1\pmod n$. The primitive divisor theorem of $\mathcal{L}_n$ ensures that for all primes $n\geq 5$, $\mathcal{L}_n$ has a primitive divisor except for finitely many pairs $(\alpha, \bar{\alpha})$ all of which are given in \cite[Table 1]{BHM01}.  
Since $d\in\{1,2,11,19,22,38,209,418\}$, using this table one sees that $\mathcal{L}_p$ has a primitive divisor except for $p=5, 7$. On the other hand, $\mathcal{L}_5$ has a primitive divisor except for $\alpha\in\lbrace\frac{1+\sqrt{-11}}{2}, 6+\sqrt{-19}\rbrace$ and
 $\mathcal{L}_7$ has a primitive divisor except for $\alpha=\frac{1+\sqrt{-19}}{2}$. These exceptions give solutions to \eqref{pe1}.

We first consider $\alpha=\frac{1+\sqrt{-11}}{2}$ for $p=5$. This gives $y=3$ and thus \eqref{pe1} reduces to \eqref{pe5}. By Lemma \ref{pl1}, the corresponding solution of \eqref{pe1} is $(x,y,k, \ell, m, p)=(1,3,1,2,0,5)$. Similarly for $\alpha=6+\sqrt{-19}$ one gets $y=55$ and hence \eqref{pe1} reduces to \eqref{pe6}. Therefore using Lemma \ref{pl2} the corresponding solution is $(x,y,k, \ell, m, p)=(22434,55,0, 0, 1,5)$.  Analogously for $\alpha=\frac{1+\sqrt{-19}}{2}$; \eqref{pe1} coincides with \eqref{pe7}
 and thus by Lemma \ref{pl3} we conclude that it does not give any solution to \eqref{pe1}. 
  
 We now assume that $\mathcal{L}_p$ has a primitive prime divisor $q$. Then $q\equiv \pm 1 \pmod p$. Applying this fact in \eqref{pe4}, we get $q = 11, 19$. 
 Since $p\geq 5$ is a prime such that  $q\equiv \pm1 \pmod p$, the only choice is  $p=5$.  Again  the sign in $q\equiv \pm 1\pmod p$ coincides with $\left( \frac{-d}{q}\right)$ and thus we have,
 $$
 d=\begin{cases} 
 1, 11& \text{ if } q=19,\\
 2, 19 & \text{ if } q=11.
 \end{cases}
 $$
We divide the remaining part of the proof in four cases.

\subsection*{Case I:  $(d, q)=(1, 19)$} 
In this case \eqref{pe3} becomes 
$$
x+z\sqrt{-1}=(u+v\sqrt{-1})^5
$$
and equating the imaginary parts gives
\begin{equation} \label{pe8}
v(5u^{4}-10u^{2}v^{2}+v^{4})=2^a11^b19^c.
\end{equation}
Thus $v\in\{\pm 1, \pm 2^{a_0}, \pm 11^{b_0}, \pm 2^{a_1}11^{b_1}\}$ for some integers $1\leq a_0, a_1\leq a$ and $1\leq b_0, b_1\leq b$ (since $\gcd(q,(\alpha-\bar{\alpha})^2)$ $= $  $ \gcd(19, -4v^2) $ $=1$).
We note that $y=u^2+dv^2=u^2+v^2$ is odd. 

Let $v=\pm 1$ and this forces $a=0$ (follows from \eqref{pe8}) and then \eqref{pe8} becomes
$$
5u^{4}-10u^{2}+1=\pm 11^b19^c.
$$
As $y$ is odd, $u$ must be  even and thus reading this equation modulo $8$ the sign on the right hand side  is positive and both $a$ and $b$ are either odd or even. Thereafter we write the above equation as, 
$$
5u^4-10u^2+1=\begin{cases}w^2 & \text{ if } b\equiv c\equiv 0\pmod 2,\\
 209 w^2 & \text{ if } b\equiv c\equiv 1\pmod 2,
\end{cases}
$$
where $w=\begin{cases} 11^{b/2}19^{c/2}& \text{ if } b\equiv c\equiv 0\pmod 2,\\
11^{(b-1)/2} 19^{(c-1)/2} & \text{if } b\equiv c\equiv 1\pmod 2.
\end{cases}$ 
\\
Using MAGMA, we see that $(u, w)=(0,1)$ is the only integer solution of $5u^4-10u^2+1=w^2$ which gives $x=0$.   

Again reading $5u^4-10u^2+1=209w^2$ modulo $5$, one gets $4w^2\equiv 1\pmod 5$ which shows that $w\equiv \pm 2\pmod 5$. This is a contradiction to $w\equiv \pm 1 \pmod 5$.  

We now consider $v=2^{a_0}$. Then \eqref{pe8} becomes 
$$
5u^{4}-10u^{2}v^2+v^4=\pm 2^{a-a_0}11^b19^c.
$$
Now $u$ is odd as $y=u^2+v^2$ is odd. Reading the last equation modulo $8$, we see that $a-a_0=0$, $b+c$ is odd and the sign on the right hand side is negative. Thus it becomes
\begin{equation}\label{pe9}
5u^{4}-10u^{2}v^{2}+v^{4}=-11^{b}19^{c}.
\end{equation}
If $b=2b_2+1$ for some non-negative integer $b_2$ then $c=2c_2$ for some positive integer $c_2$ and thus \eqref{pe9} can be written as
$$
-11V^{2}=5U^{4}-10U^{2}+1.
$$ 
Here
$$
(U,V)=\left( \frac{u}{v},\frac{11^{b_2}19^{c_2}}{v^{2}}\right) 
$$ 
and $(U,V)$ is a $\{2\}$-integral point on the above quartic curve and using MAGMA one can show that there is no such point. 

On the other hand, if $b=2b_3$ for some positive integer $b_3$ then $c=2c_3+1$ for some positive integer $c_3$ and thus \eqref{pe9} can be written as
$$
-19V^{2}=5U^{4}-10U^{2}+1,
$$
where $(U,V)=\left( \frac{u}{v},\frac{11^{b_3}19^{c_3}}{v^{2}}\right) $ is $\{2\}$-integral point on this quartic curve.  Using MAGMA, we obtain $(U,V)=\left(\pm \frac{1}{2},\pm \frac{1}{4}\right)$ which gives a previously obtained solution.

Next let $v=\pm 11^{b}$. In this case, \ref{pe8} becomes
$$5u^{4}-10u^{2}v^2+v^4=\pm 2^{a}11^{b-b_0}19^c.$$
Since $y=u^2+v^2$ is odd, $u$ is even. Reading this equation modulo $8$, we arrive at $1\equiv \pm 2^a3^{b-b_0+c}\pmod 8$ which is not possible as $a\geq 1$.

We now consider the remaining case $v=\pm 2^{a_0}11^{b_0}$.  In this situation \ref{pe8} becomes
$$
5u^{4}-10u^{2}v^{2}+v^{4}=\pm 2^{a-a_0}11^{b-b_0}19^c.
$$
Here $u$ is odd and thus reading this equation modulo $8$, we see that $a=a_0$, $b-b_0+c$ is odd and the sign on the right hand side is negative. Thus 
\begin{equation}\label{pe10}
5u^{4}-10u^{2}v^{2}+v^{4}=-11^{b-b_0}19^{c}.
\end{equation}
If $b-b_0=2b_4+1$ for some positive integer $b_4$ then $c=2c_4$ for some positive integer $c_4$ and thus \eqref{pe10} can be rewritten as
$$
-11V^{2}=5U^{4}-10U^{2}+1,$$ where 
$$
(U,V)=\begin{cases} \left( \frac{u}{v},\frac{11^{b_4+{b_0}/2}19^{c_2}}{v^{2}}\right) & \text{if } 2\mid b_0,\\
\left( \frac{u}{v},\frac{11^{b/2+b_4}19^{c_2}}{v^{2}}\right) & \text{if } 2\nmid b_0.
\end{cases}
$$ 
Clearly $(U,V)$ is  $\{2,11\}$-integral points on the above quartic curve and using MAGMA, we see that there is no such point. 

On the other hand, if $c=2c_5+1$ for some positive integer $c_5$ then $b-b_0=2b_5$ for some positive integer $b_5$ and thus \eqref{pe10} can be rewritten as
$$
-19V^{2}=5U^{4}-10U^{2}+1,
$$
where $(U,V)=\left( \frac{u}{v},\frac{11^{b_5}19^{c_5}}{v^{2}}\right) $ is $\{2,3\}$-integral point on this quartic curve and MAGMA gives $(U,V)=\left(\pm \frac{1}{2},\pm \frac{1}{4}\right)$. These lead to $(u,v)=(\pm 1, \pm 2)$ and $a=1, b=0, c=1$. All these together show that $(x,y,k,\ell,m,p)=(41, 5, 2,0,2,5)$ is a solution of \eqref{pe1}.

\subsection*{Case II:  $(d, q)=(11,19)$}
In this case, \eqref{pe3} becomes 
$$
x+z\sqrt{-11}=(u+v\sqrt{-11})^5.
$$
Now equating the imaginary parts, it gives,
\begin{equation*} 
v(5u^{4}-110u^{2}v^{2}+121v^{4})=2^a11^b19^c.
\end{equation*}
This shows that $v=\pm 2^{a_0}11^{b_0}$ for some integers $0\leq a_0 \leq a$ and $0\leq b_0\leq b$ as $\gcd(q,(\alpha-\bar{\alpha})^2)$ $= $  $ \gcd(19, -4v^2) $ $=1$. Thus 
the above equation becomes
\begin{equation}\label{pe11} 
5u^{4}-110u^{2}v^2+121v^4=\pm 2^{a-a_0}11^{b-b_0}19^c.
\end{equation}
We first assume that $v$ is even, then $u$ is odd as $y=u^2+11v^2$ is odd. Therefore reading \eqref{pe11} modulo $8$ one sees that 
$a=a_0$, $b-b_0+c$ is odd and the sign on the right hand side  is negative. Thus 
\begin{equation}\label{pe12}
5u^{4}-110u^{2}v^{2}+121v^{4}=-11^{b-b_0}19^{c}.
\end{equation}
If $b-b_0=2B_0+1$ for some non-negative integer $B_0$ then $c=2C_0$ for some positive integer $C_0$ and thus \eqref{pe12} can be read as
$$
-11V^{2}=5U^{4}-110U^{2}+121,
$$ 
where 
$$
(U,V)=\left( \frac{u}{v},\frac{11^{B_0}19^{C_0}}{v^{2}}\right)
$$ 
is a $\{2,11\}$-integral point on the above quartic curve. Using MAGMA, we conclude that there do no exist such points. 

Again if $c=2C_1+1$ for some non-negative integer $C_1$ then $b-b_0=2B_1$ for some positive integer $B_1$ and thus \eqref{pe12} can be read as
$$
-19V^{2}=5U^{4}-10U^{2}+1,
$$
where $(U,V)=\left( \frac{u}{v},\frac{11^{B_1}19^{C_1}}{v^{2}}\right) $ is $\{2,11\}$-integral point on this quartic curve. As before using MAGMA one concludes  it does not give any solution.

When $v$ is odd, that is $v=\pm 11^{b_0}$ then $u$ is even and $a_0=0$. Therefore reading \eqref{pe11} modulo $8$ pertains to
$1\equiv \pm 2^a3^{b-b_0+c}\pmod 8$. This further implies that $a=0, b-b_0+c$ is odd and sign is positive. These data reduces  \eqref{pe11} to
$$
5U^2-110U^2+121=AV^2,
$$
where $$(U,V,A)=\begin{cases}\left( \frac{u}{v}, \frac{11^{\frac{b-b_0}{2}}19^{\frac{c}{2}}}{v^2}, 1\right) & \text{if }2\mid b-b_0, 2\mid c,\\
\left( \frac{u}{v}, \frac{11^{\frac{b-b_0}{2}-1}19^{\frac{c}{2}-1}}{v^2}, 209\right) & \text{if }2\nmid b-b_0, 2\nmid c.\end{cases}$$
 We see that $(U,V)\in\{(0, \pm 11), (\pm 1, \pm 4)\}$ only when $A=1$ (using MAGMA). These values do not add to any further solution to \eqref{pe1}.

\subsection*{Case III: $(d, q)=(2, 11)$}
In this case, \eqref{pe3} becomes 
$$
x+z\sqrt{-2}=(u+v\sqrt{-2})^5
$$
and equating the imaginary parts gives,
\begin{equation*} 
v(5u^{4}-20u^{2}v^{2}+4v^{4})=2^a11^b19^c.
\end{equation*}
This shows that $v=\pm 2^{a_0}19^{c_0}$ for some integers $0\leq a_0 \leq a$ and $0\leq c_0\leq c$ as $\gcd(q,(\alpha-\bar{\alpha})^2)= \gcd(11, -4v^2) =1$. Thus 
the above equation reduces to
\begin{equation*}
5u^{4}-20u^{2}v^2+4v^4=\pm 2^{a-a_0}11^b19^{c-c_0}.
\end{equation*}
Since $y=u^2 + 2v^2$ is odd, $u$ is also odd. Thus reading the above relation  modulo $4$, we see that $a=a_0$ and  it becomes,
\begin{equation*}
5u^{4}-20u^{2}v^2+4v^4=\pm11^b19^{c-c_0}.
\end{equation*}
Let us assume that $a_0\geq 1$. Then $v$ is even and thus reading this equation modulo $8$ one gets $5\equiv \pm 3^{b+c-c_0}\pmod 8$. Similarly if $a_0=0$ then $5\equiv \pm 3^{b+c-c_0}\pmod 8$. This shows that $b+c-c_0$ is odd and sign is negative. Hence we have
\begin{equation}\label{pe13}
5u^{4}-20u^{2}v^2+4v^4=-11^b19^{c-c_0}.
\end{equation}
If $b=2b_1+1$ for some non-negative integer $b_1$ then $c-c_0=2c_1$ for some positive integer $c_1$ and thus \eqref{pe13} can be read as
$$
-11V^{2}=5U^{4}-20U^{2}+4,$$ where 
$$
(U,V)=\left( \frac{u}{v},\frac{11^{b_1}19^{c_1}}{v^{2}}\right)
$$ 
is a $\{2,19\}$-integral points on the above quartic curve. As before using MAGMA, we obtain $(U,V)\in\{(\pm 15/8, \pm 41/64), (\pm 1, \pm 1), (\pm 1/2, \pm 1/4)\}$ and thus $(u,v)\in\{(\pm 1, \pm 1), (\pm 1, \pm 2)\} $ with $b_1=c_1=0$. These further imply $(y,a,b,c)\in\{(3,0,1,0),(9,1,1,0)\}$ and thus the corresponding solutions of \eqref{pe1} are $(x,y,k,\ell,m, p)\in\{(1,3,1,2,0,5),(241,9,3,2,0,5)\}$.

Now if $c-c_0=2c_2+1$ for some non-negative integer $c_2$ then $b=2b_2$ for some positive integer $b_2$ and thus \eqref{pe13} takes the shape
$$
-19V^{2}=5U^{4}-20U^{2}+4,
$$
where $(U,V)=\left( \frac{u}{v},\frac{11^{b_2}19^{c_2}}{v^{2}}\right) $ is $\{2,19\}$-integral point on this quartic curve. One can easily conclude by using MAGMA  that there do not exist  no such points.

\subsection*{Case IV: $(d, q)=(19, 11)$} 
In this case, \eqref{pe3} becomes 
$$
x+z\sqrt{-19}=(u+v\sqrt{-19})^5
$$
and equating the imaginary parts gives,
\begin{equation*} 
v(5u^{4}-190u^{2}v^{2}+361v^{4})=2^a11^b19^c.
\end{equation*}
This shows that $v=\pm 2^{a_0}19^{c_0}$ for some integers $0\leq a_0 \leq a$ and $0\leq c_0\leq c$ as $\gcd(q,(\alpha-\bar{\alpha})^2)= \gcd(11, -4v^2) =1$. Thus 
the above equation reduces to
\begin{equation}\label{pe14}
5u^{4}-190u^{2}v^2+361v^4=\pm 2^{a-a_0}11^b19^{c-c_0}.
\end{equation}
We first assume that $a_0=0$. Then $v=19^{c_0}$ is odd and thus $u$ is even since $y=u^2+19v^2$ is odd. Therefore \eqref{pe14} becomes
$$
5u^{4}-190u^{2}v^2+361v^4=\pm 2^a11^b19^{c-c_0}.
$$
Reading this modulo $8$ one gets  $a=0$, $b+c-c_0$ is even and the sign is positive. Thus it reduces to
$$
5u^{4}-190u^{2}v^2+361v^4= 11^b19^{c-c_0}.
$$
This can be rewritten as 
$$
5U^4-190U^2+361=BV^2,
$$
where
$$
(U,V,B)=\begin{cases} \left(\frac{u}{v}, \frac{11^{b/2}19^{(c-c_0)/2}}{v^2},1\right) & \text{if } b\equiv c-c_0\equiv 0\pmod 2,\\
\left(\frac{u}{v}, \frac{11^{b/2-1}19^{(c-c_0)/2-1}}{v^2}, 209\right) & \text{if } b\equiv c-c_0\equiv 1\pmod 2.\end{cases}
$$
Using MAGMA, we get $(U,V)\in\{(0,\pm 19), (\pm 6, \pm 1)\}$ along with $b=0, c=c_0$ only when $B=1$. However $(U,V)=(\pm 6,\pm 1)$ is possible and it  gives $(u,v)=(\pm 6, \pm 1)$. This further shows that $c=0$ and $y=55$ and correspondingly  the solution to \eqref{pe1} is $(x,y,k,\ell,m,n)=(22434,55,0,0,1,5)$. 

The remaining case is when  $a_0\geq 1$. In this case, $u$ is odd and hence reading \eqref{pe14} modulo $8$ gives $a=a_0$, $b+c-c_0$ is odd and the sign is negative. Thus \eqref{pe14} reduces to
$$
5u^{4}-190u^{2}v^2+361v^4= -11^b19^{c-c_0}.
$$
This can be written as
$$
5U^4-190U^2+361=CV^2,
$$ 
where
$$
(U,V,C)=\begin{cases}(\frac{u}{v}, \frac{11^{(b-1)/2}19^{(c-c_0)/2}}{v^2}, -11)& \text{if } 2\nmid b,\\
(\frac{u}{v}, \frac{11^{b/2}19^{(c-c_0-1)/2}}{v^2}, -19)& \text{if } 2\mid b.
\end{cases}
$$
Here $(U,V)$ is $\{2,19\}$-integral points. But using MAGMA we see that there are no such points on this curve.
\end{proof}

\begin{proof}[\bf{Proof of Theorem \ref{thm}}]
The proof follows from Proposition \ref{3p}, Corollary \ref{cor3.1}, Proposition \ref{4p} and Proposition \ref{mp} except  when $n$ is a multiple of any prime $p\geq 5$. 

Let $n=pt$ for a prime $p\geq 5$ and an integer $t\geq 2$. Then \eqref{men} can be read as
\begin{equation}\label{gc}
x^2+2^k11^\ell19^m=Y^p,~ x,Y\geq 1, \gcd(x,Y)=1, k, \ell, m\geq 0, ~p\geq 5, 
\end{equation}
where $y^t=Y$.  Further applying Proposition \ref{mp} one can conclude that $Y= 3,5,9,55$. This forces $y=3$ as $y^t=Y$ with $t\geq 2$ and the corresponding solution is $(x,y,k, \ell, m, n)=(241,3,3,2,0,10)$ (see, Table \ref{Tp}). 
\end{proof}

\section{Concluding remarks} 
Here we discuss some special cases of Theorem \ref{thm} which help to deduce some known results and also give some new ones.

We first deduce the result of Cang\"{u}l et al. \cite{CDLPS10}, where they considered the Diophantine equation
\begin{equation}\label{eqd1}
x^2 + 2^k11^\ell = y^n; ~ ~ x,y\geq 1, \gcd(x,y)=1, k, \ell\geq 0.
\end{equation}
Substituting $m=0$ in Theorem \ref{thm}, we get their result \cite[Theorem 1]{CDLPS10}. Similarly, one can deduce  \cite[Corollary 2]{CDLPS10} which gives the solutions of $x^2+11^\ell=y^n$, from Theorem \ref{thm} by putting $k=m=0$. 
We note that the solutions of \eqref{eqd1} can also be obtained from \cite[Theorem 1.1]{CDILS13} by taking $b=0$.  

Soydan et al. \cite{SUZ12}
determined all the solution of 
 $$
 x^2+2^a19^b=y^n, ~ x, y\geq 1, \gcd(x,y)=1, n\geq 3. 
 $$  
 As a particular case of Theorem \ref{thm} by putting $\ell=0$ in \eqref{men}, we recover their main result. We can also deduce \cite[Corollary 2.1]{SUZ12} from Theorem \ref{thm} by taking $k=\ell=0$.

\section*{acknowledgement}
A part of this manuscript was completed while second author was visiting Professor Kotyada Srinivas at the Institute of Mathematical Sciences, Chennai. He is grateful to Professor Srinivas for his valuable comments and for stimulating environment at the Institute of Mathematical Sciences, Chennai. The authors are grateful to Harish-Chandra Research Institute (HRI) for providing facilities to carry out this work. The third author is thankful to Malaviya National Institute of Technology, Jaipur for providing facilities/financial help to carry out this work. The authors are grateful to Professor Paul Voutier and Dr. Nguyen Xuan Tho for pointing out a couple of misprints in the previous version. The authors are thankful to the anonymous referee for his/her valuable comments and suggestions which have helped improving the presentation immensely.

\end{document}